\documentclass[12pt]{amsart}
\usepackage{palatino,hhline, bbm}
\usepackage{arydshln} 
\usepackage{amsthm, amsmath}
\usepackage{amssymb} 
\usepackage{amscd}
\usepackage{mathrsfs}
\usepackage[hmargin=3cm, vmargin=3cm]{geometry}
\usepackage[breaklinks,bookmarksopen,bookmarksnumbered]{hyperref}
\usepackage[all]{xy}

\theoremstyle{plain}
\newtheorem*{thm}{Theorem}
\newtheorem{lem}{Lemma}
\newtheorem{prop}{Proposition}
\newtheorem{cor}{Corollary}

\theoremstyle{remark}

\newcommand\pr{\noindent\textit{Proof} : }
\def\lr#1{\langle {#1} \rangle}

\newcommand\im{\operatorname{Im}}
\newcommand\Ker{\operatorname{Ker}}

\newcommand\Pic{\operatorname{Pic}}

\newcommand\Z{\mathbb{Z}}

\newcommand\C{\mathbb{C}}
\renewcommand\P{\mathbb{P}}

\renewcommand\O{\mathcal{O}}

\newcommand\iso{\vbox{\hbox to .8cm{\hfill{$\scriptstyle\sim$}\hfill}
\nointerlineskip\hbox to .8cm{{\hfill$\longrightarrow $\hfill}} }}
\newcommand\abs[1]{\lvert {#1}\rvert}

\begin{document}
\title{Maximal variation of linear systems}
\author[Arnaud Beauville]{Arnaud Beauville}
\address{Universit\'e C\^ote d'Azur\\
CNRS -- Laboratoire J.-A. Dieudonn\'e\\
Parc Valrose\\
F-06108 Nice cedex 2, France}
\email{arnaud.beauville@univ-cotedazur.fr}
 
\begin{abstract}
Let $X$ be a smooth projective variety over $\C$, and $L$ a line bundle on $X$. We say that the linear system $\abs{L}$ has \emph{maximal variation} if its elements have the maximum number $\dim\abs{L}$ of moduli.
We discuss  some cases where this situation is expected: hypersurfaces, double coverings of $\P^n$, K3 and hyperk\"ahler manifolds, and abelian varieties.

\end{abstract}
\maketitle 
\section*{Introduction}

Let $X$ be a smooth projective variety over $\C$, and let $L$ be a line bundle on $X$. We assume that a general hypersurface $S$ in the linear system $\abs{L}$ is smooth and connected. We say that $L$ has \emph{maximal variation} if there are only finitely many $S'$ in $\abs{L}$ isomorphic to $S$. If there is a reasonable moduli space $\mathscr{M}$ for $S$ (this will be the case for most of our examples), this is equivalent to say that the natural rational map $\abs{L}\dasharrow \mathscr{M}$
 is generically finite, or equivalently that its image has dimension $\dim\abs{L}$. 

There is a simple differential criterion for this to happen, namely  $H^0(T_{X|S})=0$ (\S 1). However this   turns out to be surprisingly difficult to check. We will discuss some cases where maximal variation is expected: hypersurfaces (\S 2),  double coverings of $\P^n$ (\S 3),  K3 surfaces and hyperk\"ahler manifolds (\S 4), and abelian varieties (\S 5). We show indeed that the variation is maximal  in the last three cases, while we have only partial results in the two first cases. Finally we discuss some more examples and counter-examples (\S 6 and 7).

After the publication of this paper on arXiv, Bricalli and Pirola have obtained (by a completely different method) a very general result \cite{B-P}:
\begin{thm}\label{BP}
Let $X$ be a smooth projective variety of dimension $n$, with $h^{n,0}\neq 0$, $h^{n-1,0}=0$. Any very ample line bundle on $X$ has maximal variation.
\end{thm}
As we will indicate below, this improves in particular some of our results.

\section{A differential criterion}

\begin{prop}\label{gen}
$(i)$  If $H^0(S, T_{X|S})=0$, $L$ has maximal variation.

$(ii)$ If $H^1(X,\O_X)=0$ and $H^0(S,T_S)=0$, the converse holds.
\end{prop}

\pr Let $S$ be a general element of $\abs{L}$, and let $k:\mathscr{S}\rightarrow B$ be a Kuranishi family for $S$. There is an open neighborhood $U$ of $[S]$ in $\abs{L}$ and a map $\mu :U\rightarrow B$ such that the pull back of $k$ by $\mu $ is the family of hypersurfaces parameterized by $U$. $L$ has maximal variation if and only if $\mu $ is a local immersion at $[S]$, that is, the tangent map $T\mu $ is injective.

The exact sequence $0\rightarrow T_S\rightarrow T_{X|S}\rightarrow L_{|S}\rightarrow 0$ gives a cohomology exact sequence
\[0\rightarrow H^0(T_S)\rightarrow H^0(T_{X|S})\rightarrow H^0(L_{|S}) \xrightarrow{\ \partial \ }H^1(T_S) \,.  \]
Let $s$ be a section of $L$ with zero locus  $S$. The exact sequence $0\rightarrow \O_X\xrightarrow{\ s\ }L\rightarrow L_{|S}\rightarrow 0$ gives a cohomology exact sequence
\[0\rightarrow H^0(L)/(\C\!\cdot\! s) \xrightarrow{\ r \ } H^0(L_{|S})\rightarrow H^1(\O_X)\,.\]The tangent space to $\abs{L}$ at $S$ is $H^0(L)/(\C\!\cdot \!s)$, and the tangent space to $B$ at $\mu (S)$ is $H^1(T_S)$; by {\cite[(12.4)]{K-S}, the tangent map to $\mu $ at $[S]$ is the composition \[T\mu :H^0(L)/(\C\!\cdot \!s) \xrightarrow{\ r \ } H^0(L_{|S})\xrightarrow{\ \partial \ }H^1(T_S)\,.\]
If $H^0(T_{X|S})=0$, $\partial  $ is injective; since $r $ is injective, $T\mu $ is injective. Conversely if $T\mu $ is injective, under the hypotheses of $(ii)$ 
$r $ is an isomorphism, so $\partial $ is injective and $H^0(T_{X|S})=0$ since $H^0(T_S)=0$.\qed

\medskip	

\section{Hypersurfaces}
Let $X\subset \P^n$ be a smooth hypersurface of degree $d\geq 3$, defined by an equation $F=0$. We assume $n\geq 3$, and $d\geq 4$ if $n=3$. The \emph{Jacobian ring} of $X$ is the graded artinian ring $R:=\C[x_0,\ldots ,x_n]/(F'_0,\ldots ,F'_n)$, where we write $F'_i=\dfrac{\partial F}{\partial x_i} \,\cdot $
\begin{prop}\label{hyper}
$\O_X(1)$ has maximal variation if and only if the multiplication map  $\times \ell: R_{d-1}\rightarrow R_d$ is injective for $\ell$ general in $R_1$.
\end{prop}

The property that  $\times \ell: R_{p-1}\rightarrow R_p$ is of maximal rank for $\ell$ general in $R_1$ and for all $p$ is called the \emph{weak Lefschetz property};  it is conjectured that it holds for all   rings of the form $\C[x_0,\ldots ,x_n]/(P_0,\ldots ,P_n)$, with $P_0,\ldots ,P_n$ homogeneous with no common zeroes in $\P^n$ (see for instance \cite{J-M}). This would imply that  $\O_X(1)$ has always maximal variation. We have only a weaker result:

\begin{cor}\label{cor}
$\O_X(1)$ has   maximal variation in the following cases:

$1)\ d\geq n+2$;

$2)\ X$ is a surface in $\P^3\  ($with $d\geq 4)$;

$3)\ X$ is a cubic threefold;

$4)\ X$ is general \emph{(among hypersurfaces of degree $d$ in $\P^n$).}
\end{cor}

\pr 1) is proved in \cite{BMR}, improving by 1 the bound in \cite{PRT}. This implies 2) for $d\geq 5$; for $d=4$ and $n=3$, the weak Lefschetz property is proved in \cite{BMMN}. 

3) is proved in \cite{B}, as a consequence of \cite{A-R}.

4) The Fermat hypersurface $\sum X_i^d=0$ in $\P^n$ satisfies the weak Lefschetz property (in all degrees), because its Jacobian ring is isomorphic to the cohomology ring $H^{*}((\P^{d-2})^{n},\C)$. Therefore the same holds for $X$ in a Zariski open subset of $\abs{\O_{\P^n}(d)}$.
\qed

\medskip	
\noindent \emph{Remark}$.-$ 1) is improved into  $d\geq n+1$ in \cite{B-P}, as a consequence of the theorem mentioned in the introduction.
This  theorem applies also to all non-Fano (smooth) complete intersections.

\medskip	
\noindent\emph{Proof of the Proposition} : We essentially repeat the proof of \cite{B}, with one more detail for the case $n=3$. 
 
Let $S$ be a smooth hyperplane section of $X$; we choose the coordinates so that it is defined by $x_0=0$.
From the exact sequence
\[0\rightarrow T_{X|S}\rightarrow T_{\P^n|S}\rightarrow \O_S(d)\rightarrow 0\]
we see that $H^0(T_{X|S})$ is the kernel of the homomorphism $\varphi : H^0(T_{\P^n|S})\rightarrow H^0(\O_S(d))$. 

From the commutative diagram of Euler exact sequences
\[\xymatrix{0 \ar[r] & H^0(\O_{\P^n}) \ar[r]\ar[d]^{\wr} & H^0(\O_{\P^n}(1)^{n+1})\ar[r]\ar@{->>}[d]& H^0(T_{\P^n})\ar[r]\ar[d]&0&&\\
0 \ar[r] & H^0(\O_{S}) \ar[r] &H^0 (\O_{S}(1)^{n+1})\ar[r]& H^0(T_{\P^n|S})\ar[r]&H^1(\O_S)\ar[r]^<<<<{u}&H^1(\O_S(1)^{n+1})
}\]
we get an exact sequence
\[H^0(\P^n, T_{\P^n}) \xrightarrow{\ r \ } H^0(S,T_{\P^n|S}) \rightarrow  H^1(S,\O_S) \xrightarrow{\ u \ } H^1(S,\O_S(1)^{n+1})\,.\]

If $n\geq 4$ we have $H^1(S,\O_S)=0$; if $n=3$ (so that $S$ is a curve), the homomorphism $u$ is the transpose of the map $H^0(S, \O_S(d-4))^4 \xrightarrow{\ (x_0,\ldots ,x_3)\ } H^0(S, \O_S(d-3))$, which is surjective since $d\geq 4$; thus $u$ is injective.  In each case we conclude that  the restriction map $r $ is surjective. 
It follows that $H^0(S,T_{\P^n|S})$ is generated by the elements $x_i\,\dfrac{\partial }{\partial x_j}$ for $i\geq 1$, $j\geq 0$, with the Euler relation $\sum  \limits_{i\geq 1}x_i \,\dfrac{\partial }{\partial x_i}=0$. 

Assume that $\times x_0:R_{d-1}\rightarrow R_d$ is injective. Let $V\in H^0(S,T_{\P^n|S})$; we can write $V=\sum\limits_{i=0}^n L_i  \,\dfrac{\partial }{\partial x_i}$, where the $L_i$ are linear forms in $x_1,\ldots ,x_n$.
The homomorphism $\varphi : H^0(T_{\P^n|S})\rightarrow H^0(\O_S(d))$ maps $V$ to $\sum L_i F'_i$. 
 If $\varphi (V)=0$ , the form $\sum L_iF'_i$ vanishes on $S$, hence there exists  $G\in H^0(\P^n,\O_{\P}(d-1))$ and $a\in\C$ such that
$\sum L_i F'_i= x_0 G+aF$. Since $dF=\sum x_iF'_i$, it follows that $x_0G$ is zero in $R_d$. By hypothesis this implies $G=0$ in $R_{d-1}$, that is, $G=\sum a_iF'_i$ for some constants $a_i$. Therefore \[\sum_{i=0}^n (L_i-a_ix_0- bx_i)F'_i= 0\qquad \mbox{with }b:=a/d\,.\]

Since $X$ is smooth, $(F'_0,\ldots ,F'_n)$ is a regular sequence in $\C[x_0,\ldots ,x_n]$, hence
\allowbreak$L_i-a_ix_0- bx_i=0$ for each $i$. Since $x_0$ does not appear in $L_i$, we get :

\smallskip	
\centerline{$\bullet$ For $i\geq 1$, $a_i=0$ and $L_i=bx_i$;\qquad
$\bullet$ for $i=0$, $L_0=0$.}
\smallskip	

Hence $V= b\sum \limits_{i\geq 1}x_i \,\dfrac{\partial }{\partial x_i}=0 $. Thus $\varphi $ is injective, therefore $H^0(T_{X|S})=0$, and $\O_X(1)$ has maximal variation by Proposition \ref{gen} $(i)$.

\medskip	
Conversely, suppose that $\O_X(1)$ has maximal variation. We have $H^1(X,\O_X)=0$ and $H^0(S, T_S)=0$ because of our hypotheses on $n$ and $d$,  
so $H^0(S,T_{X|S})=0$ by Proposition \ref{gen}$(ii)$, and $\varphi : H^0(T_{\P^n|S})\rightarrow H^0(\O_S(d))$ is injective. 
Let $G \in H^0(\P^n,\O_{\P}(d-1))$ such that $x_0G=0$ in $R_d$; there exist linear forms $L_i$ such that $x_0G=\sum L_iF'_i$. Replacing $G$ by $G-\sum a_iF'_i$, where $a_i$ is the coefficient of $x_0$ in $L_i$, we can assume that the $L_i$ are linear forms in $x_1,\ldots ,x_n$. Then the vector field $V= \sum L_i\,\dfrac{\partial }{\partial x_i}$ satisfies $\varphi (V)=0$, hence $V=0$. This means that there exists $c\in\C$ such that $L_i=cx_i$ for $i>0$ and $L_0=0$. Then $x_0G= c\sum\limits_{i>0}x_iF'_i= c(dF-x_0F'_0)$, hence
$cdF= x_0(G+cF'_0)$. Since $F$ is irreducible, this implies $c=0$, hence $G=0$. \qed

\medskip	
We now consider the linear system $\O_X(e)$ for $e\geq 2$. Recall that $\dim X\geq 2$, $d:=\deg X\geq 3$, and $d\geq 4$ if $\dim X=2$.
\begin{prop}
$\O_X(e)$ has maximal variation in the following cases :

$1)\  \O_X(1)$  has maximal variation;

$2)\ e\geq d$. 
\end{prop}
(In particular, each of the conditions 1) to 3) of Corollary \ref{cor} implies that $\O_X(e)$ has maximal variation for all $e\geq 1$).

\smallskip	
\pr The proof of Proposition \ref{hyper} works identically for  $\O_X(e)$, replacing $x_0$ by a general form $h$ in $H^0(\O_{\P^n}(e))$.
We find that $\O_X(e)$ has maximal variation if $\times h: R_{d-e}\rightarrow R_d$ is injective for $h$ general in $R_e$.

1) Let $\ell$ be a general element of $R_1$. If $\O_X(1)$  has maximal variation, $\times \ell: R_{d-1}\rightarrow R_d$ is injective (Proposition \ref{hyper}). Then $\times \ell: R_{p-1}\rightarrow R_p$ is injective for $p\leq d$ (if $\ell x=0$ for $x\in R_{p-1}$, we have for every $y$ in $R_{d-p}$ $\ell (xy)=0$ , hence $xy=0$; thus
$x\cdot R_{d-p}=0$, which  implies $x=0$ since $R$ is Gorenstein). Therefore $\times \ell^{e}: R_{d-e}\rightarrow R_d$ is injective, and so is $\times h$ for $h$ general in $R_e$. 

2) If $e>d$, we have $R_{d-e}=0$ so there is nothing to prove. If $e=d$, 
we just observe that $R_d\neq 0$, since
the socle of $R$ is in degree $ (n+1)(d-2)\geq d$.\qed

\medskip	
\section{Double coverings of $\P^n$}

Let $\pi :X\rightarrow \P^n$ a double covering, branched along a smooth hypersurface $B\subset \P^n$ of (even) degree $d$. We assume $n\geq 2$, $d\geq 4$, and $d\geq 6$ if $n=2$. Let $F=0$ be an equation of $B$, and
let $R=\C[x_0,\ldots ,x_n]/(F'_0,\ldots ,F'_n)$ be the Jacobian ring of $F$.

\begin{prop}\label{double}
The line bundle $\pi ^*\O_{\P}(1)$ has maximal variation if and only if the multiplication map $\times \ell: R_{d-1}\rightarrow R_d$ is injective for $\ell$ general in $R_1$.
\end{prop}

Using   \cite{B-P}, we get:
\begin{cor}\label{cordouble}
If $d\geq n+1\ ($and $d\geq 6$ if $n=2)$, $\pi ^*\O_{\P}(1)$ has maximal variation.

\end{cor}

\noindent\emph{Proof of the Proposition} : The proof is inspired by  \cite[Proposition 5.2]{D-H}, which treats the case $n=2$, $d=6$.

Put $\delta :=d/2$, and $\P=\P^n$. We have $\pi _*\Omega ^1_X= \Omega ^1_{\P}\oplus \Omega ^1_{\P}(\log B)(-\delta )$ \cite[Lemma 3.16 (d)]{E-V}. Using Grothendieck duality this gives 
\[\pi _*T_X = T_{\P}(-\delta )\oplus T_{\P}(-\log B)\]
where $T_{\P} (-\log B)$ is the sheaf of vector fields on $\P$ tangent to $B$, dual to $\Omega ^1_{\P}(\log B)$.

\smallskip	
Let $H$ be a general hyperplane in $\P$, and $S:=\pi ^*H$.
We have \[\pi_*T_{X|S}=(\pi _*T_X)_{|H}= T_{\P}(-\delta )_{|H}\oplus T_{\P}(-\log B)_{|H}\,,\]hence the vanishing of $H^0(T_{X|S})$ is equivalent to $H^0( T_{\P}(-\delta )_{|H})=H^0(T_{\P}(-\log B)_{|H})=0$.
\smallskip	

1) Consider the Euler exact sequence restricted to $H$:
\[ 0\rightarrow \O_{H}(-\delta )\rightarrow \O_{H}(-\delta +1)^{n+1}\rightarrow T_{\P}(-\delta )_{|H}\rightarrow 0\,.\]
If $n\geq 3$, we have $H^1(\O_{H}(-\delta ))=0$, hence $H^0(T_{\P}(-\delta )_{|H})=0$. If $n=2$ (so $H\cong \P^1$), the homomorphism $H^1(\O_{H}(-\delta ))\rightarrow H^1(\O_{H}(-\delta +1)^3)$ is the transpose of the map $H^0(\O_H(\delta -3)^3)\xrightarrow{\  (x_0,x_1,x_2)\ }H^0(\O_H(\delta -2)) $, which is surjective since $\delta \geq 3$; hence again $H^0(T_{\P}(-\delta )_{|H})=0$.

\smallskip	
2) Consider the commutative diagram of exact sequences 
\[\xymatrix{0\ar[r] & \O_{\P}\ar[r]\ar@{=}[d] &\O_{\P}(1)^{n+1}\ar[r]\ar[d]^{\varphi }& T_{\P}
\ar[r]\ar[d]^{\psi}& 0\\
0\ar[r] & \O_{\P}\ar[r]^{dF} & \O_{\P}(d) \ar[r] & N_{B/\P}=\O_B(d)\ar[r]&0
}\]where $\varphi (L_0,\ldots ,L_n)=\sum L_iF'_i$. Since $\Ker \psi=T_{\P}(-\log B)$, the snake lemma gives an exact sequence
\[0\rightarrow T_{\P}(-\log B)\rightarrow \O_{\P}(1)^{n+1} \xrightarrow{\ \varphi \ } \O_{\P}(d)\rightarrow 0\,.\]
 Therefore we have $H^0( T_{\P}(-\log B)_{|H})=0$ if and only if $\varphi _{|H}: H^0(\O_H(1)^{n+1})\rightarrow H^0(\O_H(d))$ is injective. Choose the coordinates so that $H$ is defined by $x_0=0$.
Let  $L_0,\ldots ,L_n$ be linear forms in $x_1,\ldots ,x_n$ such that $\varphi _{|H}(L_0,\ldots ,L_n)=0$. This means that 
 $\sum L_iF'_i=x_0G$ for some $G\in H^0(\O_{\P}(d-1))$. If $\times x_0: R_{d-1}\rightarrow R_d$ is injective, 
we get $G=\sum a_iF'_i$ for some scalars $a_i$, thus $\sum (L_i-a_ix_0)F'_i=0$. As before this implies 
$L_i=a_ix_0$, and since  $x_0$ does not appear in $L_i$, $L_i=0$.

Conversely, suppose that $\varphi _{|H}$ is injective. Let $G\in R_{d-1}$ such that $x_0G=0$ in $R_d$: there exist linear forms $L_0,\ldots ,L_n$ such that $x_0G=\sum L_iF'_i$. 
This implies that the restrictions $\bar{L}_i$ of $L_i$ to $H$ satisfy $\varphi _{|H}(\bar{L}_0,\ldots ,\bar{L}_n)=0$, hence $\bar{L}_i=0$ for all $i$. Therefore $L_i=b_ix_0$ for some scalars $b_i$, so that $G=\sum b_iF'_i$ is zero in $R_{d-1}$.\qed

\smallskip	
\noindent\emph{Remarks}$.-$ 1) Propositions \ref{double} and \ref{hyper} have the following  curious consequence. The hypersurface $B$ embeds into $X$; under the hypotheses of Proposition \ref{double},
\emph{$\pi ^*\O_{\P}(1)$ has maximal variation if and only if its restriction to $B$ has maximal variation}. I don't know a direct proof of this equivalence. 

2) As in \S 2, the proof extends to the case of the line bundle $\pi ^*\O_{\P}(e)$ for $e>1$: we get that $\pi ^*\O_{\P}(e)$ \emph{has maximal variation when $\pi ^*\O_{\P}(1)$ does, or when} $e\geq d$.

\medskip	
\section{K3 and hyperk\"ahler}

\begin{prop}
Any ample line bundle on a K3 surface has maximal variation.
\end{prop}

This is \cite[Corollary 2]{Ba} (which uses \cite[Proposition 5.2]{D-H}).\qed

\medskip	
The case of higher-dimensional hyperk\"ahler manifolds turns out to be   easier:
\begin{prop} Let $X$ be a hyperk\"ahler manifold of dimension $>2$. Any ample line bundle $L$ on $X$  has maximal variation.
\end{prop}
\pr Let $S\in\abs{L}$. The exact sequence $0\rightarrow T_X\otimes L^{-1}\rightarrow T_X\rightarrow T_{X|S}\rightarrow 0$ gives rise to a cohomology exact sequence 
$H^0(T_X)\rightarrow H^0(T_{X|S})\rightarrow H^1(T_X\otimes L^{-1})$. Since $T_X\cong \Omega ^1_X$, we have $H^0(T_X)=0$ and $H^1(T_X\otimes L^{-1})=0$ by the Akizuki-Nakano vanishing theorem, hence $H^0(T_{X|S})=0$.\qed

\smallskip	
When $L$ is very ample, the two above results  follow also from \cite{B-P}.
\medskip	
\section{Abelian varieties}

\begin{prop}\label{ab}
Let $X $ be an abelian variety. Any   ample line bundle $L$ on $X$ has maximal variation.
\end{prop}

\smallskip	
\pr As in Proposition \ref{gen}, we need to prove that the map $T\mu :H^0(L)/(\C\cdot s) \rightarrow H^1(T_S)$ is injective.
Let $s$ be a section of $L$ with zero locus  $S$. The exact sequences \[0\rightarrow T_S\rightarrow T_{X|S}\rightarrow L_{|S}\rightarrow 0 \qquad\mbox{and}\qquad 0\rightarrow \O_X \xrightarrow{\ s\ } L\rightarrow L_{|S}\rightarrow 0\] give rise to a diagram
\[\xymatrix{ & H^0(L)/(\C\!\cdot\! s)\ar[d]^r\ar[dr]^{T\mu}&\\
H^0(T_{X|S})\ar[r]^\varphi \ar[rd]^u &H^0(L_{|S})\ar[r]\ar[d]^{\partial }& H^1(T_S)\\
& H^1(\O_X)&}\]where the horizontal and vertical lines are exact.

We claim that $u$ is an isomorphism; this implies that $\im \varphi \cap \im r=(0)$, hence that $T\mu $ is injective. Since the vector bundle $T_X$ is trivial, it is equivalent to consider the composition $\tilde{u}: H^0(T_X)\iso H^0(T_{X|S})\xrightarrow{\ u\ }H^1(\O_X)$.

\begin{lem}
$\tilde{u}: H^0(T_X)\rightarrow H^1(\O_X)$   is the cup-product with the Chern class $c_1(L)\in H^1(\Omega ^1_X)$.
\end{lem}
\pr We can find a finite covering $(U_{\alpha })_{\alpha \in I}$ of $X$ such that $L_{|U_{\alpha }}$ is free, generated by a section $e_{\alpha }$; on $U_{\alpha }\cap U_{\beta }$ we have $e_{\alpha }=g_{\alpha \beta }e_{\beta }$ for an invertible function $g_{\alpha \beta }$.
We write $s_{|U_{\alpha }}=s_{\alpha }e_{\alpha }$ for a function $s_{\alpha }$ on $U_{\alpha }$; then $s_{\beta }=g_{\alpha \beta }s_{\alpha }$. 

Let $V\in H^0(T_X)$.  
Then  $\varphi (V)$ is given on $U_{\alpha }\cap S$ by $(V\cdot s_{\alpha })e_{\alpha }$. 
This expression lifts to $U_{\alpha }$, and $u(V)=\partial \varphi (V)$ is the class of the \v{C}ech 1-cocycle $(\alpha ,\beta )\mapsto c_{\alpha \beta }$ such that :
 \[ c_{\alpha \beta }\,s=  (V\cdot s_{\beta })e_{\beta }-(V\cdot s_{\alpha })e_{\alpha }= (V\cdot g_{\alpha \beta })g_{\alpha \beta }^{-1}\,s\ \mbox{ in }\, U_{\alpha }\cap U_{\beta }\, ,\ \mbox{ hence } \ c_{\alpha \beta }=\lr{V,g_{\alpha \beta }^{-1}\,dg_{\alpha \beta }}\,.\]The class of the cocycle $(\alpha ,\beta )\mapsto g_{\alpha \beta }^{-1}\,dg_{\alpha \beta }$ in $H^1(X,\Omega ^1_X)$ is the Chern class $c_1(L)$, hence the lemma.\qed

\begin{lem}
$\tilde{u}$  is an isomorphism.
\end{lem} 

\pr  We put $g=\dim X$, and choose an isomorphism $\Omega ^g_X\iso \O_X$. This gives an isomorphism  $T_X\iso \Omega ^{g-1}_X$; for $\alpha \in H^0(\Omega ^{1}_X)$, the map $\lr{-,\alpha }:T_X\rightarrow \O_X$ is identified with the 
 exterior product $ \wedge \alpha :\Omega ^{g-1}_X\rightarrow \Omega ^{g}_X$. Therefore $\tilde{u}$ is identified with the cup-product $\cup\, c_1(L): H^0(\Omega_X ^{g-1})\rightarrow H^1(\Omega ^g_X)$. By the Lefschetz theorem the cup-product $\cup\, c_1(L):H^{g-1}(X,\C)\rightarrow H^{g+1}(X,\C)$ is an isomorphism, and it induces an isomorphism from $H^{g-1,0}$ onto $H^{g,1}$, hence our assertion.\qed

\medskip	
\section{Some more examples}

We collect here some other examples of line bundles $L$ with maximal variation. We always assume that the general element of  $\abs{L}$ is smooth and connected.

$\bullet$ If $H^0(X,T_X)=0$ and $L$ is ample, $L^m$ has maximal variation for $m\gg 0$. Indeed for any $S\in\abs{L^m}$ there is an exact sequence $H^0(T_X)\rightarrow H^0(T_{X|S})\rightarrow H^1(T_X\otimes L^{-m})$, and the right hand term vanishes for  $m\gg 0$ since $L$ is ample.

$\bullet$ If $\Omega ^1_X$ is ample, every line bundle $L$ on $X$ has maximal variation. Indeed for any $S\in \abs{L}$ the quotient $\Omega ^1_{X|S}$ is ample, hence does not admit a nonzero map to $\O_S$. 

\smallskip	
For surfaces we have a more precise result:
\begin{prop}
Let $X$ be a  surface  of general type such that $\Omega ^1_X$ is generated by its global sections. Then any
very ample line bundle $L$ on $X$ has maximal variation.
\end{prop}
\pr Since $\Omega ^1_X$ is globally generated, the zero locus $Z(\omega )$ of a general section $\omega \in H^0(\Omega ^1_X)$    is finite and reduced (see e.g. \cite[Remark 6]{K}), and nonempty since $c_2>0$. Let $p\in Z(\omega )$. Since $L$ is very ample, we can find a smooth curve   $C\in\abs{L}$ such that $C\cap Z(\omega )=\{p\} $.  We have an exact sequence
\[0\rightarrow K^{-1}_{X|C}(p)\rightarrow T_{X|C}\xrightarrow{\ \omega \ }\O_C(-p)\rightarrow 0\,;\]
it suffices to prove $\deg (K^{-1}_{X|C}(p))<0$, that is, $K_X\cdot C>1$.

We  observe that $X$  is a minimal surface: if $R$ is a smooth  curve in $X$, $\Omega ^1_R$ is a quotient of $\Omega ^1_X$, hence it must be globally generated, which implies $g(R)\geq 1$. 
Then by the index theorem we have $(K_X\cdot C)^2\geq K_X^2 \cdot C^2>1$, hence $K_X\cdot C>1$ as required.\qed

\begin{prop}
Let $X$ be a surface of general type such that 
$\Pic(X)=\Z\cdot [K_X]$ and that the Chern classes of $X$ satisfy 
$c_1^2>c_2$.   The line bundles $K_X^{r}$, for $r\geq 1$,  have maximal variation.
\end{prop}

\pr The proof is borrowed from \cite{B1}. Let $C$ be a smooth curve in $\abs{K_X^{r}}$. 
Suppose that $T_{X|C}$ has a nonzero section; this gives a surjective homomorphism $\Omega ^1_X\rightarrow \O_C(-E)$, where $E$ is an effective divisor. The kernel $F$ of this homomorphism 
 is a rank 2 vector bundle with \[c_1(F)=(r-1)c_1\ ,\quad  c_2(F)= c_2-rc_1^2-\deg(E)\, ,\]
 \[ c_1^{2}(F)-4c_2(F)\geq (r+1)^2c_1^2-4c_2\geq 4(c_1^2-c_2)>0\,.\leqno{\mbox{hence}}\]
By \cite[\S 10]{Bo} $F$ is unstable: there is an exact sequence
\[0\rightarrow K_X^{a}\rightarrow F\rightarrow \mathscr{I}_ZK_X^{b}\rightarrow 0\]with $a\geq b$, and $Z$  a finite subscheme of $X$. 

This gives $c_1(F)=-(a+b)c_1$, $c_2(F)=\deg(Z)+abc_1^2$, hence
 \[c_1^{2}(F)-4c_2(F)=(a-b)^2c_1^2-4\deg(Z)\,.\] 
 Comparing with the previous expression  gives 
\[(a-b)^2c_1^2 \geq  c_1(F)^2-4c_2(F)\geq   (r+1)^2c_1^2-4c_2> (r^2+2r-3)c_1^2\geq (r-1)^2c_1^2 \]
Since $a\geq b$, we conclude $a-b> r-1$, hence $a\geq 1$. This gives a nonzero homomorphism $K_X\rightarrow \Omega ^1_X$, hence a nonzero section of $T_X$, which is impossible. \qed

\begin{prop}
Let $X=X_1\times X_2$, where $X_1$ and $X_2$ are smooth projective varieties of dimension $\geq 2$, and $H^0(X_i,T_{X_i})=0$.
 Then any  ample line bundle on $X$ has maximal variation.
\end{prop} 
\pr
 Let $p_i$ be the projection of $X$ onto $X_i$.
The restriction $\bar{p}_i:S\rightarrow X_i$ of $p_i$ to an ample divisor $S\subset X$ is surjective with connected fibers: 
indeed for $x\in X_i$,  $\bar{p}_i^{-1}(x)=S\cap (\{x_i\}\times  X_{1-i}) $ is an ample divisor in  $\{x_i\}\times  X_{1-i}$, hence nonempty and connected. Therefore we have $(\bar{p}_i)_{*}\O_S=\O_{X_i}$, and \[H^0((p_i^*T_{X_i})_{|S})=H^0(\bar{p}_i^*T_{X_i})=H^0(T_{X_i})=0 \ ,\ \mbox{ hence }\ H^0(T_{X|S})=0\,. \qed\]

\medskip	
\section{Some counter-examples}

$\bullet$ Suppose that a connected  linear group $G$ acts non trivially on $X$.  Then $G$ acts trivially on $\Pic(X)$, hence acts on each linear system $\abs{L}$ on $X$. If this action is nontrivial, $L$ does not have maximal variation. This is the case for instance if $L$ is very ample, since one can find an element of $\abs{L}$ passing through a point of $X$ and not containing its orbit.

\smallskip	
$\bullet$ Let $C$ be a smooth projective curve,  $E$ a rank 2 very ample vector bundle on $C$ (that is, $\O_{\P(E)}(1)$ is very ample), $X=\P(E)$, $L=\O_{\P(E)}(1)$.
A general element  of $\abs{L}$ is isomorphic to $C$,  so $\abs{L}$ has no moduli.

\smallskip	
$\bullet$ There are varieties of general type admitting a  base point free pencil with no moduli. Take for instance the Fermat hypersurface $F: \sum\limits_{i=0}^n x_i^d=0$ in $\P^n$. Consider the rational map $p:F\dasharrow \P^1$ given by $p(x_0,\ldots ,x_n)= \allowbreak (x_0,x_{1})$. Let $\hat{F}\rightarrow F$ be the blowing up of the codimension $2$ subvariety of $F$ given by $x_0=x_{1}=0$.  Then $p$ extends to a morphism $\hat{p}:\hat{F}\rightarrow \P^1$ whose general fiber is isomorphic to the Fermat hypersurface of degree $d$ in $\P^{n-1}$. The pencil  $\abs{\hat{p}^*\O_{\P^1}(1)}$ has  no moduli.

\medskip	
These examples lead to the following :
\smallskip	

\noindent\textbf{Conjecture}$.-$ \emph{If $X$ is not uniruled,  every ample line bundle  on $X$ has maximal variation}.

\end{document}